\documentclass[graybox]{svmult}
\usepackage{mathptmx,helvet,courier,graphicx}
\usepackage{psfrag,euscript,amsmath,amsfonts,amssymb,amscd,amsopn}
\usepackage[ansinew]{inputenc}

\newcommand{\be}{\beta}

\newcommand{\vfi}{\varphi}
 \newcommand{\RR}{{\mathbb R}}

  \def \cS {\mathcal{S}}
\def \cF {{\mathcal F}}

\newcommand{\cO}{{\EuScript O}}

\newcommand{\EE}{{\mathbb E}}

\newcommand{\NN}{{\mathbb N}}

\newcommand{\dist}{\operatorname{dist}}

\newcommand{\Leb}{\operatorname{Leb}}

\newcommand{\sgn}{\operatorname{sgn}}

\newcommand{\suporte}{\operatorname{supp}}

\begin{document}

\title{Lorenz-like chaotic attractors revised}

\author{V{\'{\i}}tor Ara\'ujo and Maria Jos\'e Pacifico}
\institute{V\'{\i}tor Ara\'ujo \at Instituto de Matem\'a\-tica,
Universidade Federal do Rio de Janeiro,\\
C. P. 68.530, 21.945-970, Rio de Janeiro, RJ-Brazil
\\\emph{and} \\
Centro de Matem\'atica da
  Universidade do Porto, \\Rua do Campo Alegre 687, 4169-007
  Porto, Portugal
\\
\email: {vitor.araujo@im.ufrj.br and vdaraujo@fc.up.pt}
\and
Maria Jos\'e Pacifico \at Instituto de Matem\'atica,
Universidade Federal do Rio de Janeiro,\\
C. P. 68.530, 21.945-970 Rio de Janeiro, Brazil
\\
\email: {pacifico@im.ufrj.br and pacifico@impa.br}}

\maketitle

\abstract{ We describe some recent results on the dynamics
  of singular-hyperbolic (Lorenz-like) attractors $\Lambda$
  introduced in \cite{MPP04}: (1) there exists an invariant
  foliation whose leaves are forward contracted by the flow;
  (2) there exists a positive Lyapunov exponent at every
  orbit; (3) attractors in this class are expansive and so
  sensitive with respect to initial data; (4) they have zero
  volume if the flow is $C^2$, or else the flow is globally
  hyperbolic; (5) there is a unique physical measure whose
  support is the whole attractor and which is the
  equilibrium state with respect to the center-unstable
  Jacobian; (6) the hitting time associated to a geometric
  Lorenz attractor satisfies a logarithm law; (7) the rate
  of large deviations for the physical measure on the
  ergodic basin of a geometric Lorenz attractor is
  exponential.}



\section{Introduction}
In this note $M$ is a compact boundaryless $3$-manifold and
${\cal X}^1(M)$ denotes the set of $C^1$ vector fields on
$M$ endowed with the $C^1$ topology. Moreover $\Leb$ denotes
\emph{volume} or \emph{Lebesgue measure}: a normalized
volume form given by some Riemannian metric on $M$.  We also
denote by $\dist$ the Riemannian distance on $M$.

The notion of singular hyperbolicity was introduced
in~\cite{MPP98,MPP04} where it was proved that any $C^1$ robustly
transitive set for a $3$-flow is either a singular
hyperbolic attractor or repeller.

A compact invariant set
$\Lambda$ of a $3$-flow $X \in {\cal X}^1(M)$ is
an \emph{attractor} if there exists a neighborhood $U$ of
$\Lambda$ (its isolating neighborhood) such that
\begin{align*}
  \Lambda=\bigcap_{t>0} X^t(U)
\end{align*}
and there exists $x\in\Lambda$ such that $X(x)\neq\vec0$ and
whose positive orbit $\{X^t(x):t>0\}$ is dense in $\Lambda$.

We say that a compact invariant subset is \emph{singular
  hyperbolic} if all the singularities in $\Lambda$ are
hyperbolic, and the tangent bundle $T\Lambda$ decomposes in
two complementary $DX^t$-invariant bundles $E^s \oplus
E^{cu}$, where: $E^s$ is one-dimensional and uniformly
contracted by $DX^t$; $E^{cu}$ is bidimensional, contains the
flow direction, $DX^t$ expands area along $E^{cu}$ and
$DX^t\mid E^{cu}$ dominates $DX^t\mid E^s$ (i.e.  any eventual
contraction in $E^s$ is stronger than any possible
contraction in $E^{cu}$), for all $t>0$.  

We note that the presence of an equilibrium  together with
regular orbits accumulating on it prevents any invariant set
from being uniformly hyperbolic, see
e.g.~\cite{BR75}. Indeed, in our $3$-dimensional setting a
compact invariant subset $\Lambda$ is uniformly hyperbolic
if the tangent bundle $T\Lambda$ decomposes in \emph{three}
complementary $DX^t$-invariant bundles $E^s \oplus E^{X}
\oplus E^u$, each one-dimensional, $E^X$ is the flow
direction, $E^s$ is uniformly contracted and $E^u$ uniformly
expanded by $DX^t$, $t>0$. This implies the continuity of
the splitting and the presence of a non-isolated
equilibrium point in $\Lambda$ leads to a discontinuity in
the splitting dimensions.

In the study of the asymptotic behavior of orbits of a flow
$X \in {\cal X}^1(M)$, a fundamental problem is to
understand how the behavior of the tangent map $DX$
determines the dynamics of the flow $X_t$.  The main
achievement along this line is the uniform hyperbolic
theory: we have a complete description of the dynamics
assuming that the tangent map has a uniformly hyperbolic
structure since~\cite{BR75}.

In the same vein, under the assumption of singular
hyperbolicity, one can show that at each point there exists a
strong stable manifold and that the whole set is foliated by
leaves that are contracted by forward iteration. In
particular this shows that any robust transitive attractor
with singularities displays similar properties to those of
the geometrical Lorenz model.  It is also possible to show
the existence of local central manifolds tangent to the
central unstable direction.  Although these central
manifolds do not behave as unstable ones, in the sense that
points on them are not necessarily asymptotic in the past,
the expansion of volume along the central unstable
two-dimensional direction enables us to deduce some
remarkable properties.

We shall list some of these properties that give us a nice
description of the dynamics of a singular hyperbolic
attractor.

\section{The geometric Lorenz attractor}
\label{sec:geoLorenz}

Here we briefly recall the construction of the geometric
Lorenz attractor \cite{ACS, GW79}, that is the more representative
example of a singular-hyperbolic attractor.

In 1963 the meteorologist Edward Lorenz published in the
Journal of Atmospheric Sciences \cite{Lo63} an example of a
parametrized polynomial system of differential equations
\begin{align}
  \label{e-Lorenz-system}
\dot x &= a(y - x)&
\quad
&a = 10 \nonumber
\\
\dot y &= rx -y -xz&
\quad
&r =28
\\
\dot z &= xy - bz&
\quad
&  b = 8/3 \nonumber
\end{align}
as a very simplified model for thermal fluid convection,
motivated by an attempt to understand the foundations of
weather forecast.

The origin $\sigma=(0,0,0)$ is an equilibrium of saddle type for the vector 
field defined by equations (\ref{e-Lorenz-system}) with real eigenvalues 
$\lambda_i$, $i\leq 3$ satisfying
\begin{equation}
\label{eigenvalues}
  \lambda_2<\lambda_3<0<-\lambda_3<\lambda_1.
\end{equation}
(in this case $\lambda_1\approx 11.83$ , $\lambda_2\approx
-22.83$, $\lambda_3=-8/3$).

\begin{figure}[htpb]
  \sidecaption
  \includegraphics[scale=0.7]{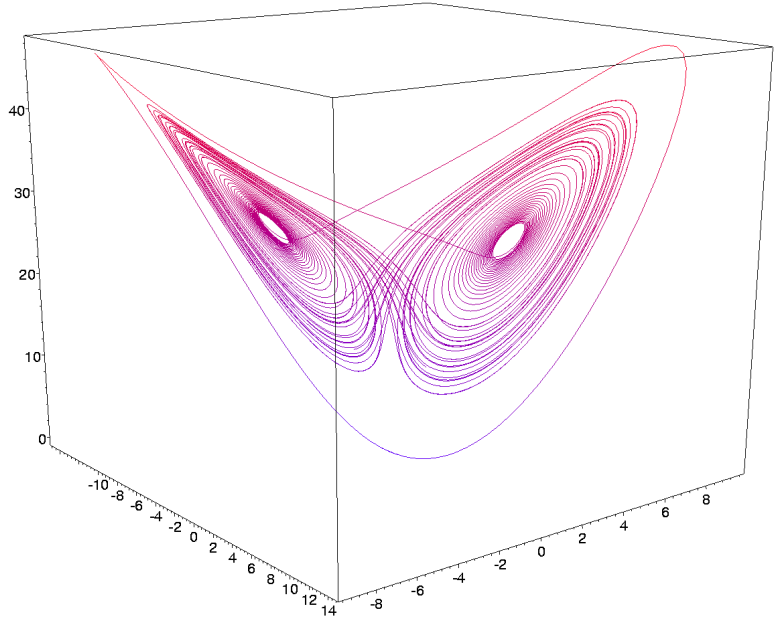}
  \caption{A view of the Lorenz attractor calculated numerically}
  \label{fig:view-lorenz-attract}
\end{figure}

Numerical simulations performed by Lorenz for an open
neighborhood of the chosen parameters suggested that almost
all points in phase space tend to a {\em chaotic attractor},
whose well known picture is presented in
Figure~\ref{fig:view-lorenz-attract}.  The \emph{chaotic
  feature} is the fact that trajectories converging to the
attractor are {\em sensitive with respect to initial data}:
trajectories of any two nearby points are driven apart under
time evolution.

Lorenz's equations proved to be very resistant to
rigorous mathematical analysis, and also presented
serious difficulties to rigorous numerical study.
Indeed, these two main difficulties are:
\begin{description}
\item[\emph{conceptual}:] the presence of an equilibrium point
  at the origin accumulated by regular orbits of the flow
  prevents this atractor from being hyperbolic
  \cite{AraPac07},
\item[\emph{numerical}:] the presence of an equilibrium
  point at the origin, implying that solutions slow down as
  they pass near the origin, which means unbounded return
  times and, thus, unbounded integration errors.
\end{description}
Moreover the attractor is \emph{robust}, that is, the
features of the limit set persist for all nearby vector
fields. More precisely, if $U$ is an isolating neighborhood
of the attractor $\Lambda$ for a vector field $X$, then
$\Lambda$ is \emph{robustly transitive} if, for all vector
fields $Y$ which are $C^1$ close to $X$, the corresponding
$Y$-invariant set
\begin{align*}
  \Lambda_Y(U)=\bigcup_{t>0} Y^t(U)
\end{align*}
also admits a dense positive $Y$-orbit.  We remark that
\emph{the persistence of transitivity}, that is, the fact
that, for all nearby vector fields, the corresponding limit
set is transitive, \emph{implies a dynamical
  characterization of the attractor}, as we shall see.

These difficulties led, in the seventies, to the
construction of a geometric flow presenting a similar
behavior as the one generated by equations
(\ref{e-Lorenz-system}). Nowadays this model is known
as {\em geometric Lorenz flow}.  Next we briefly describe
this construction, see \cite{ACS, GW79} for full details.

We start by observing that under some non-resonance
conditions, by the results of Sternberg~\cite{St58}, in a
neighborhood of the origin, which we assume to contain the
cube $[-1,1]^3 \subset \RR^3$, the Lorenz equations are
equivalent to the linear system $(\dot x, \dot y, \dot
z)=(\lambda_1 x,\lambda_2 y, \lambda_3 z)$ through
smooth conjugation, thus
\begin{align}\label{eq:LinearLorenz}
X^t(x_0,y_0,z_0)=
(x_0e^{\lambda_1t}, y_0e^{\lambda_2t}, z_0e^{\lambda_3t}),
\end{align}
where $\lambda_1\approx 11.83$ , $\lambda_2\approx -22.83$,
$\lambda_3=-8/3$ and $(x_0, y_0, z_0)\in\RR^3$ is an
arbitrary initial point near $(0,0,0)$.

Consider $S=\big\{ (x,y,1) : |x|\le
{\scriptstyle{1/2}},\quad |y|\le{\scriptstyle{1/2}}\big\}$ and
\begin{align*}
S^-&=\big\{ (x,y,1)\in S : x<0 \big\},&
\qquad
S^+&=\big\{ (x,y,1)\in S : x>0 \big\}\quad\text{and}
\\
S^*&=S^-\cup S^+=S\setminus\ell, & \text{where}\quad
\ell &=\big\{(x,y,1)\in S : x=0 \big\}.
\end{align*}
Assume that $S$ is a global transverse section to the flow so that
every trajectory eventually crosses $S$ in the direction of
the negative $z$ axis.

Consider also $\Sigma=\{ (x,y,z) : |x|=1
\}={\Sigma}^-\cup{\Sigma}^+$ with ${\Sigma}^{\pm}=\{ (x,y,z)
: x=\pm 1\}$. 

For each $(x_0,y_0,1)\in S^*$ the time $\tau$
such that $X^{\tau}(x_0,y_0,1)\in\Sigma$ is given by
$$\tau(x_0)=-\frac{1}{\lambda_1}\log{|x_0|},$$ which depends
on $x_0\in S^*$ only and is such that $\tau(x_0)\to+\infty$
when $x_0\to 0$. This is one of the reasons many standard
numerical algorithms were unsuited to tackle the Lorenz
system of equations. Hence we get
(where $\sgn(x)=x/|x|$ for $x\neq0$)
\begin{equation}\label{L} 
X^\tau(x_0,y_0,1)=
\big( \sgn(x_0),  y_0e^{\lambda_2\tau},
e^{\lambda_3\tau}\big)
=
\big( \sgn(x_0),
y_0|x_0|^{-\frac{\lambda_2}{\lambda_1}},
|x_0|^{-\frac{\lambda_3}{\lambda_1}}\big).
\end{equation}
Since $0<-\lambda_3<\lambda_1<-\lambda_2$, we  have
$0<\alpha=-\frac{\lambda_3}{\lambda_1} <1
<\beta=-\frac{\lambda_2}{\lambda_1}$.
Let $L:S^*\to\Sigma$ be such that $ L(x,y)=\big(
y|x|^\beta,|x|^\alpha \big)$ with the convention that
$L(x,y)\in{\Sigma}^+$ if $x>0$ and $L(x,y)\in{\Sigma}^-$ if $x<0$.
\begin{figure}[h]
\sidecaption
\includegraphics[width=7.5cm]{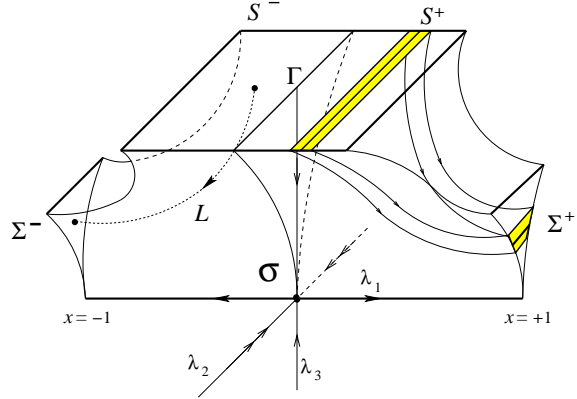}
\caption{\label{L3Dcusp}Behavior near the origin.}
\end{figure}
It is easy to see that $L(S^\pm)$ has the shape of a
triangle without the vertex $(\pm 1,0,0)$.
In fact the vertex $(\pm 1,0,0)$ are
cusp points at the boundary of each of these sets.
The fact that $0<\alpha<1<\beta$ together with equation (\ref{L}) imply
that $L(\Sigma^\pm)$ are uniformly compressed in the $y$-direction.

From now on we denote by $\Sigma^\pm$ the closure of
$L(S^\pm)$.  Clearly each line segment $S^*\cap\{x=x_0\}$ is
taken to another line segment $\Sigma\cap\{z=z_0\}$ as
sketched in Figure~\ref{L3Dcusp}.

The sets $\Sigma^\pm$ should return to the cross section $S$
through a composition of a translation $T$, an expansion $E$
only along the $x$-direction and a rotation $R$ around
$W^s(\sigma_1)$ and $W^s(\sigma_2)$, where $\sigma_i$ are
saddle-type singularities of $X^t$ that are outside the cube
$[-1,1]^3$, see \cite{AraPac07}. We assume that this
composition takes line segments $\Sigma\cap\{z=z_0\}$ into
line segments $S\cap\{x=x_1\}$ as sketched in
Figure~\ref{L3Dcusp}.  The composition $T\circ E\circ R$ of
linear maps describes a vector field $V$ in a region outside
$[-1,1]^3$.  The geometric Lorenz flow $X^t$ is then defined
in the following way: for each $t\in\RR$ and each point $x
\in S$, the orbit $X^t(x)$ will start following the linear
field until $\tilde\Sigma^\pm$ and then it will follow $V$
coming back to $S$ and so on. Let us write ${\cal
  B}=\{ X^t(x), x\in S, t\in \RR^+\} $ the set where this
flow acts.  The geometric Lorenz flow is then the pair
$({\cal B}, X^t )$ defined in this way.  The set
$$
\Lambda=\cap_{t\ge 0}X^t(S)
$$
is the {\em geometric Lorenz attractor}.

\begin{figure}[htb]
  \centering
  \psfrag{e}{$\ell$}\psfrag{y}{$f(x)$}\psfrag{S}{$S$}
  \psfrag{w}{$-1$}\psfrag{z}{$0$}\psfrag{u}{$1$}
  \psfrag{l1}{$\lambda_1$}\psfrag{l2}{$\lambda_2$}\psfrag{l3}{$\lambda_3$}
  \includegraphics[width=10cm]{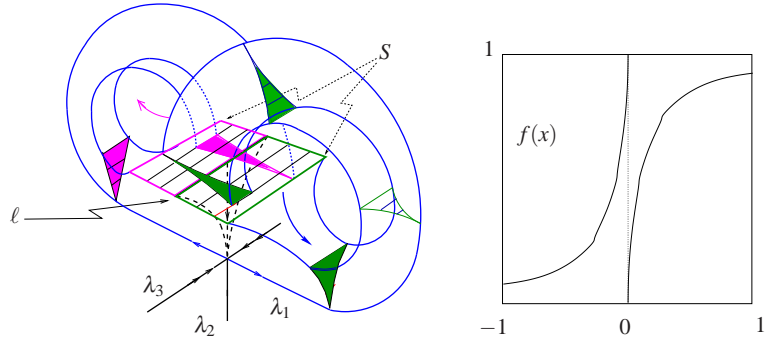}
\caption{\label{fig:GeoLorenz} The global cross-section for the geometric Lorenz
  flow and the associated 1d quotient map, the Lorenz transformation.}
\end{figure}

We remark that the existence of a chaotic attractor for the original Lorenz
system was established by Tucker with the help of
a computer aided proof (see~\cite{Tu2}).

The combined effects of $T\circ E\circ R$ and the linear flow given by
equation (\ref{L}) on lines implies that
the foliation $\cF^s$ of $S$ given by the lines $S\cap\{x=x_0\}$ is
invariant under the first return map $F:S\to S$. In another words, we have

\bigskip
$(\star)$ \/{\em for any given
leaf $\gamma$ of $\cF^s$, its image $F(\gamma)$ is
contained in a leaf of $\cF^s$}.
\bigskip

\begin{figure}[htbp]
\sidecaption
\includegraphics[width=5cm]{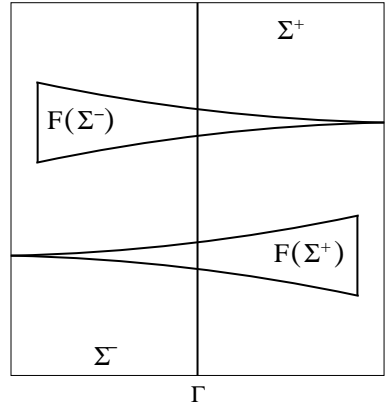}
\caption{\label{L2D}{The image $F(S^*)$.}}
\end{figure}

The main features of the geometric Lorenz flow and its first return map
can be seen at figures \ref{fig:GeoLorenz} and \ref{L2D}.

\begin{figure}[htbp]
\sidecaption
\includegraphics[width=5cm]{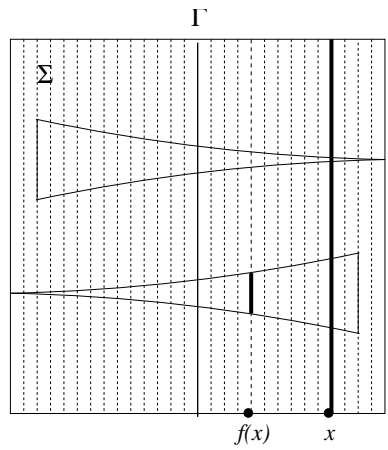}
\caption{\label{folhea}{Projection on $I$.}}
\end{figure}

The one-dimensional map $f$ is obtained quotienting over the
leaves of the stable foliation $\cF^s$ defined before.

For a detailed construction of a geometric Lorenz flow
see~\cite{AraPac07,SZ09}.

As mentioned above, a geometric Lorenz attractor is the most
representative example of a singular-hyperbolic
attractor~\cite{MPP98}.

\section{The dynamical results}
\label{sec:dynamic-results}

The study of robust attractors is inspired by the Lorenz
flow example.  Next we list the main dynamical properties of
a robust attractor.

\subsection{Robustness and singular-hyperbolicity}
\label{sec:robustn-singul-hyper}

Inspired by the Lorenz flow example we define an equilibrium
$\sigma$ of a flow $X^t$ to be {\em Lorenz-like} if the
eigenvalues $\lambda_1, \lambda_2, \lambda_3$ of
$DX(\sigma)$ are real and satisfy the relation at
(\ref{eigenvalues}):
\begin{equation*}
  \lambda_2<\lambda_3<0<-\lambda_3<\lambda_1.
\end{equation*}

These are the equilibria contained in robust attractors
naturally, since they are the only kind of equilibria in a
$3$-flow which cannot be perturbed into saddle-connections
which generate sinks or sources when unfolded.

\begin{theorem}
\label{tlorenzlike}
Let $\Lambda$ be a robustly transitive set of $X \in {\cal
  X}^1(M)$.  Then, either for $Y=X$ or $Y=-X$, every
singularity $\sigma\in\Lambda$ is Lorenz-like for $Y$ and
satisfies
$
W_Y^{ss}(\sigma)\cap\Lambda=\{\sigma\}.
$
\end{theorem}

The fact that a robust attractor does not admit sinks or
sources for all nearby vector fields in its isolating
neighborhood has several other strong consequences (whose
study was pioneered by R. Ma\~n\'e in its path to solve the
Stability Conjecture in~\cite{Man82}) which enable us to
show the following, see \cite{MPP04}.

\begin{theorem}
  A robustly transitive set for $X\in{\cal X}^1(M)$ is a
  singular-hyperbolic attractor for $X$ or for $-X$.
\end{theorem}

The following shows in particular that the notion of
singular hyperbolicity is an extension of the notion of
uniform hyperbolicity.

\begin{theorem}
\label{hypinduced}
Let $\Lambda$ be a singular hyperbolic compact set of $X\in
{\cal X}^1(M)$.  Then any invariant compact set
$\Gamma\subset \Lambda$ without singularities is uniformly
hyperbolic.
\end{theorem}

A consequence of Theorem \ref{hypinduced} is that every
periodic orbit of a singular hyperbolic set is
hyperbolic. The existence of a periodic orbit in every
singular-hyperbolic attractor was proved recently in~\cite{BM04}
and also a more general result was obtained in~\cite{AP}.

\begin{theorem}
\label{ccc}
Every singular hyperbolic attractor $\Lambda$ has a dense
subset of periodic orbits.
\end{theorem}

In the same work~\cite{AP} it was announced that every singular
hyperbolic attractor is the homoclinic class associated to
one of its periodic orbits.  Recall that the
\emph{homoclinic class} of a periodic orbit $\cO$ for $X$ is
the closure of the set of transversal intersection points of
it stable and unstable manifold:
$H(\cO)=\overline{W^u(\cO)\pitchfork W^s(\cO)}$. This result
is well known for the elementary dynamical pieces of
uniformly hyperbolic attractors. Moreover, in particular,
the geometric Lorenz attractor is a homoclinic class as
proved in~\cite{B04}.

\subsection{Singular-hyperbolicity and chaotic behavior}
\label{sec:singul-hyperb-chaoti}

Using the area expansion along the bidimensional central
direction, which contains the direction of the flow, one can
show

\begin{theorem}
  Every orbit in any singular-hyperbolic attractor has a
  direction of exponential divergence from nearby orbits
  (positive Lyapunov exponent).
\end{theorem}

Denote by $S(\RR)$ the set of surjective increasing continuous
real functions $h:\RR\to \RR$ endowed with the $C^0$
topology.  The flow $X_t$ is {\em expansive} on an invarian
compact set $\Lambda$ if for every
$\epsilon >0$ there is $\delta >0$ such that if for some $h
\in S(\RR)$ and $x,y\in\Lambda$
$$
\dist(X_t(x),X_{h(t)}y)\leq \delta
\quad\text{for all}\quad t\in\RR,
$$
then $X_{h(t_0)}(y) \in X_{[t_0-\epsilon,t_0+\epsilon]}(x)$,
for some $t_0 \in \RR$. A stronger notion of expansiveness
was introduced by Bowen-Ruelle~\cite{BR75} for uniformly
hyperbolic attractors, but equilibria in expansive sets
under this strong notion must be isolated, see
e.g. \cite{AraPac07}.

Komuro proved in~\cite{Km84} that a geometrical Lorenz attractor
$\Lambda$ is expansive.  In particular, this implies that
this kind of attractor is \emph{sensitive with respect to
  initial data}, i.e., there is $\delta>0$ such that for
any pair of distinct points $x, y \in \Lambda$, if
$\dist(X_t(x),X_t(y))<\delta$ for all $t \in \RR$, then $x$
is in the orbit of $y$. In~\cite{APPV} this was fully
extended to the singular-hyperbolic setting.

\begin{theorem}
\label{sci}
Let $\Lambda$ be a singular hyperbolic attractor of $X\in
{\cal X}^1(M)$.  Then $\Lambda$ is expansive.
\end{theorem}

\begin{corollary}
\label{sensi}
Singular hyperbolic attractors are sensitive
with respect to initial data.
\end{corollary}


%
%



\subsection{Singular-hyperbolicity, positive volume and global hyperbolicity}
\label{sec:singul-hyperb-positi}

Recently a generalization of the results of
Bowen-Ruelle~\cite{BR75} was obtained in~\cite{AAPP} showing
that a uniformly hyperbolic transitive subset of saddle-type
for a $C^{1+\alpha}$ flow has zero volume, for any
$\alpha>0$. We denote the family of all flows whose
differentiability class is at least H\"older-$C^1$ by
$C^{1+}$.

\begin{theorem}
  \label{mthm:C+attractor}
  A $C^{1+}$ singular-hyperbolic attractor has zero volume.
\end{theorem}

This can be extended to the following dichotomy.  Recall
that a \emph{transitive Anosov vector field} $X$ is a vector
field without singularities such that the entire manifold
$M$ is a uniformly hyperbolic set of saddle-type.

\begin{theorem}
  \label{mthm:volposanosov}
  Let $\Lambda$ be a singular hyperbolic attractor for a
  $C^{1+}$ $3$-dimensional vector field $X$.  Then either
  $\Lambda$ has zero volume or $X$ is a transitive Anosov
  vector field.
\end{theorem}

\section{The Ergodic Theory of singular-hyperbolic attractors}
\label{sec:ergodic-theory-singu}

The ergodic theory of singular-hyperbolic attractors is
incomplete.  Most results still are proved only in the
particular case of the geometric Lorenz flow, which
automatically extends to the original Lorenz flow after the
work of Tucker \cite{Tu2}, but demand an extra effort to
encompass the full singular-hyperbolic setting.

\subsection{Existence of a physical measure}
\label{sec:existence-physic-mea}

Another main result obtained in~\cite{APPV} is that typical
orbits in the basin of every singular-hyperbolic attractor,
for a $C^2$ flow $X$ on a $3$-manifold, have well-defined
statistical behavior, i.e. for Lebesgue almost every point
the forward Birkhoff time average converges, and it is given
by a certain physical probability measure. It was also
obtained that this measure admits absolutely continuous
conditional measures along the center-unstable directions on
the attractor. As a consequence, it is a $u$-Gibbs state and
an equilibrium state for the flow.

\begin{theorem}\label{srb}
  A $C^2$ singular-hyperbolic attractor $\Lambda$ admits a
  unique ergodic physical hyperbolic invariant probability
  measure $\mu$ whose basin covers Lebesgue almost every
  point of a full neighborhood of $\Lambda$.
\end{theorem}
Recall that an invariant probability measure $\mu$ for a
flow $X$ is \emph{physical} (or \emph{SRB}) if its \emph{basin}
\begin{align*}
  B(\mu)=
  \Big\{x\in M: 
  \lim_{T\to\infty}\frac1T\int_0^T\psi(X_t(x))\,dt=\int\psi\,d\mu,
  \forall \psi\in C^0(M,\RR)
  \Big\}
\end{align*}
has positive volume in $M$. 

Here hyperbolicity means \emph{non-uniform hyperbolicity} of
the probability measure $\mu$:
the tangent bundle over $\Lambda$ splits into a sum $T_z M =
E^s_z\oplus E^X_z\oplus F_z$ of three one-dimensional
invariant subspaces defined for $\mu$-a.e.  $z\in \Lambda$
and depending measurably on the base point $z$, where $\mu$
is the physical measure in the statement of
Theorem~\ref{srb}, $E^X_z$ is the flow direction (with zero
Lyapunov exponent) and $F_z$ is the direction with positive
Lyapunov exponent, that is, for every non-zero vector $v\in
F_z$ we have
\[
\lim_{t\to+\infty}\frac1t\log\|DX_t(z)\cdot v\|>0.
\]
We note that the invariance of the splitting implies that
$E^{cu}_z=E^X_z\oplus F_z$ whenever $F_z$ is defined.

Theorem~\ref{srb} is another statement of sensitiveness,
this time applying to the whole essentially open set
$B(\Lambda)$.  Indeed, since non-zero Lyapunov exponents
express that the orbits of infinitesimally close-by points
tend to move apart from each other, this theorem means that
most orbits in the basin of attraction separate under
forward iteration.  See Kifer~\cite{Ki88}, and
Metzger~\cite{mtz001}, and references therein, for previous
results about invariant measures and stochastic stability of
the geometric Lorenz models.

The $u$-Gibbs property of $\mu$ is stated as follows. 

\begin{theorem}
\label{thm:srbmesmo}
Let $\Lambda$ be a singular-hyperbolic attractor for a $C^2$
three-dimen\-sional flow.  Then the physical measure $\mu$
supported in $\Lambda$ has a disintegration into absolutely
continuous conditional measures $\mu_\gamma$ along
center-unstable surfaces $\gamma$ such that
$\frac{d\mu_\gamma}{dm_\gamma}$ is uniformly bounded from
above.  Moreover $\suporte(\mu)=\Lambda\,$.
\end{theorem}

Here the existence of unstable manifolds is guaranteed by
the hyperbolicity of the physical measure: the
strong-unstable manifolds $W^{uu}(z)$ are the ``integral
manifolds'' in the direction of the one-dimensional
sub-bundle $F$, tangent to $F_z$ at almost every
$z\in\Lambda$. The sets $W^{uu}(z)$ are embedded
sub-manifolds in a neighborhood of $z$ which, in general,
depend only measurably (including its size) on the base
point $z\in\Lambda$.
The \emph{strong-unstable manifold} is defined by
\[
W^{uu}(z)=\{ y\in M:
\lim_{ t\to -\infty} \dist(X_t(y),X_t(z))=0 \}
\]
and exists for almost every $z\in\Lambda$ with respect to
the physical and hyperbolic measure obtained in
Theorem~\ref{srb}. We remark that since $\Lambda$ is an
attracting set, then $W^{uu}(z)\subset\Lambda$ whenever
defined. The central unstable surfaces mentioned in the
statement of Theorem~\ref{thm:srbmesmo} are just small 
strong-unstable manifolds carried by the flow, which are
tangent to the central-unstable direction $E^{cu}$.

The absolute continuity property along the center-unstable
sub-bundle given by Theorem~\ref{thm:srbmesmo} ensures that
\[
h_\mu(X^1)=\int \log\big| \det (DX^1\mid E^{cu})  \big| \, d\mu,
\]
by the characterization of probability measures satisfying
the Entropy Formula, obtained in~\cite{LY85}.  The above
integral is the sum of the positive Lyapunov exponents along
the sub-bundle $E^{cu}$ by Oseledets
Theorem~\cite{Man87,Wa82}.  Since in the direction $E^{cu}$
there is only one positive Lyapunov exponent along the
one-dimensional direction $F_z$, $\mu$-a.e. $z$, the
ergodicity of $\mu$ then shows that the following is true.

\begin{corollary}
    \label{cor:uGibbs}
  If $\Lambda$ is a singular-hyperbolic attractor for a
  $C^2$ three-dimen\-sional flow $X^t$, then the physical
  measure $\mu$ supported in $\Lambda$ satisfies the Entropy
  Formula
\[
h_\mu(X^1)=\int\log\| DX^1\mid F_z\|\, d\mu(z).
\]
\end{corollary}

Again by the characterization of measures satisfying the
Entropy Formula, we get that
\emph{$\mu$ has absolutely continuous disintegration along
  the strong-unstable direction}, along which the Lyapunov
exponent is positive, thus \emph{$\mu$ is a $u$-Gibbs
  state}~\cite{PS82}.  This also shows that \emph{$\mu$ is
  an equilibrium state for the potential} $-\log\| DX^1\mid
F_z\|$ with respect to the diffeomorphism $X^1$. We note
that the entropy $h_\mu(X^1)$ of $X^1$ is the entropy of the
flow $X^t$ with respect to the measure $\mu$ \cite{Wa82}.

Hence we are able to extend most of the basic results on the
ergodic theory of hyperbolic attractors to the setting of
singular-hyperbolic attractors.

\subsection{Hitting and recurrence time versus local dimension
  for geometric Lorenz flows}
\label{sec:hitting}

Given $x\in M$, let $B_r(x)=\{y\in M; d(x,y)\leq r\}$ be the ball centered 
at $x$ with radius $r$.
The {\em local dimension} of $\mu$ at $x\in M$ is defined by
$$
d_\mu(x)=\lim_{r\to\infty}\frac{\log\mu(B_r(x))}{\log r}
$$
if this limit exists.
In this case $\mu(B_r(x)) \sim r^{d_\mu(x)}$.

This notion characterizes the local geometric structure of
an invariant measure with respect to the metric in the phase
space of the system, see \cite{Yo82} and \cite{Pe97}.

The existence of the local dimension for a Borel probability
measure $\mu$ on $M$ implies the crucial fact that virtually
{\em all} the known characteristics of dimension type of the
measure coincide.  The common value is a fundamental
characterisc of the fractal structure of $\mu$, see
\cite{Pe97}.

%

%

Let $x_{0}\in \RR^3$ and
\begin{equation*}
\tau _{r}^{X^t}(x,x_{0})=\inf \{t\geq 0 \,\,|\,\, X^{t}(x)\in B_r(x_{0})\}
\end{equation*}
be the time needed for the $X$-orbit of a
point $x$ to enter for the {\em first time} in a ball
$B_{r}(x_{0})$. The number $\tau_{r}^{X^t}(x,x_{0})$
is the {\em{hitting time associated to}} the flow $X^t$ and $B_r(x_0)$.
If the orbit $X^{t}$ starts at $x_{0}$ itself and we consider the
second entrance time in the ball
\begin{equation*}\label{eq-2}
\tau _{r}^{\prime }(x_{0})=\inf \{t\in \mathbb{R}^{+}:X^{t}(x_0)\in
B_{r}(x_{0}),\exists i, s.t. X^{i}(x_0)\notin B_{r}(x_{0})\}
\end{equation*}%
we have a quantitative recurrence indicator, and the number
$\tau _{r}^{\prime }(x_{0})$ is called the {\em recurrence
  time associated to} the flow $X^t$ and $B_r(x_0)$.

Now let $X^t$ be a geometric Lorenz flow, and $\mu$ its
$X^t$-invariant SBR measure.

The main result in \cite{SZ09} establishes the following

\begin{theorem}
\label{th:SZ09}
For $\mu$-almost every $x$,
{\bf{$${\bf{\lim_{r\to 0}\frac{\log \tau_r(x,x_0)}{-\log r}=d_\mu(x_0)-1.}}$$}}
\end{theorem}

Observe that the result above indicates once more the
chaoticity of a Lorenz-like attractor: it shows that
asymptotically, such attractors behave as an i.d. system.

We can always define the {\em upper} and the {\em lower}
local dimension at $x$ as
$$
{d}^+_\mu(x)=\lim \sup_{r\to\infty}\frac{\log\mu(B_r(x))}{\log r}\,,\quad\quad
{d}^-_\mu(x)=\lim \inf_{r\to\infty}\frac{\log\mu(B_r(x))}{\log r}.
$$
If $d^+(x)=d^-(x)=d$ almost everywhere the system is called
{\em exact dimensional}.  In this case many properties of
dimension of a measure coincide.  In particular, $d$ is
equal to the dimension of the measure $\mu$: $d=\inf\{\dim_H
Z; \mu(Z)=1\}.$ This happens in a large class of systems,
for example, in $C^{2}$ diffeomorphisms having non zero
Lyapunov exponents almost everywhere, \cite{Pe97}.

Using a general result proved in \cite{S06} it is also
proved in \cite{SZ09} a quantitative recurrence bound for
the Lorenz geometric flow:

\begin{theorem}
 \label{co:XSaussol}
For a geometric Lorenz flow  it holds
$$
\liminf_{r\to 0}\frac{\log \tau'_r(x)}{-\log r}={d}^-_{\mu}-1,\quad\quad\quad
\limsup_{r\to 0}\frac{\log \tau'_r(x)}{-\log r}={d}^+_{\mu}-1,
\quad \mu-a. e.\,.
$$
where $\tau'$ is the recurrence time for the flow, as defined above.
\end{theorem}

The proof of Theorem \ref{th:SZ09} is based on the following
results, proved in \cite{SZ09}.

Let $F:S \to S$ be the first return map to $S$, a global
cross section to $X^t$ through $W^s(p)$, $p$ the singularity
at the origin, as indicated at Figure \ref{fig:GeoLorenz}.
It follows that $F$ has a physical measure $\mu_F$, see
e.g.~\cite{Vi97}.  Recall that we say the system $(S, F,
\mu_F)$ has exponential decay of correlation for Lipschitz
observables if there are constants $C>0$ and $\lambda>0$,
depending only on the system such that for each $n$ it holds
\begin{equation*}
  \left|\int ~g(F^{n}(x))f(x)d\mu-\int g(x)d\mu \int f(x)d\mu \right|\leq C
  \cdot e^{-\lambda n}  \label{l1bv}
\end{equation*}%

for any Lipschitz observable $g$ and $f$ with bounded variation,
\begin{theorem}
  \label{th:1} Let $\mu_F$ an invariant physical measure for
  $F$. The system $(S, F, \mu_F)$ has exponential decay of
  correlations with respect to Lipschitz observables.
\end{theorem}

We remark that a sub-exponential bound for the decay of
correlation for a two dimensional Lorenz like map was given
in \cite{Bu83} and \cite{ACS}.

\begin{theorem} $\mu_F$ is exact, that is, $d_{\mu_F}(x)$
  exist almost every $x\in S$.
\end{theorem}

Let $x_0\in S$ and $\tau_r^S(x,x_0)$ be the time needed to
$\cO_x$ enter for the first time in $B_r(x_0)\cap S=B_{r,S}.$

\begin{theorem} $\lim_{r\to
    0}\frac{\log{\tau_{r}(x,x_0)}}{-\log{r}} = \lim_{r\to
    0}\frac{\log{\tau_{r}^S(x,x_0)}}{-\log{r}}=d_{\mu_F}(x_0).$
\end{theorem}

From the fact that the attractor is a suspension of the
support of $\mu_F$ we easily deduce the following.

\begin{theorem}
$d_\mu(x)=d_{\mu_F}(x) + 1.$
\end{theorem}

We remark that the results in this section can be extended to 
a more general class of flows described in \cite{SZ09}.
The interested reader can find the detailed proofs in this
article.

\subsection{Large Deviations for the physical measure on a
  geometric Lorenz flow}
\label{sec:large-deviat}

Having shown that physical probability measures exist, it is
natural to consider the rate of convergence of the time
averages to the space average, measured by the volume of the
subset of points whose time averages stay away from the
space average by a prescribed amount up to some evolution
time.  We extend part of the results on large deviation
rates of Kifer ~\cite{kifer90} from the uniformly hyperbolic
setting to semiflows over non-uniformly expanding base
dynamics and unbounded roof function. These special flows
model non-uniformly hyperbolic flows like the Lorenz flow,
exhibiting equilibria accumulated by regular orbits.

\subsubsection{Suspension semiflows}
\label{sec:suspens-semifl}

We first present these flows and then state the main
assumptions related to the modelling of the geometric Lorenz
attractor.

Given a H\"older-$C^1$ local diffeomorphism
$f:M\setminus\cS\to M$ outside a volume zero
non-flat\footnote{$f$ behaves like a power of the distance
  to \( \cS \): $\| Df(x) \| \approx \dist(x,\cS)^{-\be}$
  for some $\beta>0$ (see
  Alves-Arajo~\cite{alves-araujo2004} for a precise
  statement).  } singular set $\cS$, let $X^t:M_r\to M_r$ be
a \emph{semiflow with roof function $r:M\setminus\cS\to\RR$
  over the base transformation $f$}, as follows.

\begin{figure}[ht]
  \sidecaption
  \psfrag{8}{$\infty$}
  \psfrag{W}{$X^t_r(x_0,y_0)$}
  \psfrag{r1}{$-r(x_1)$}
  \psfrag{r2}{$-r(x_2)$}
  \includegraphics[width=10cm,height=4cm]{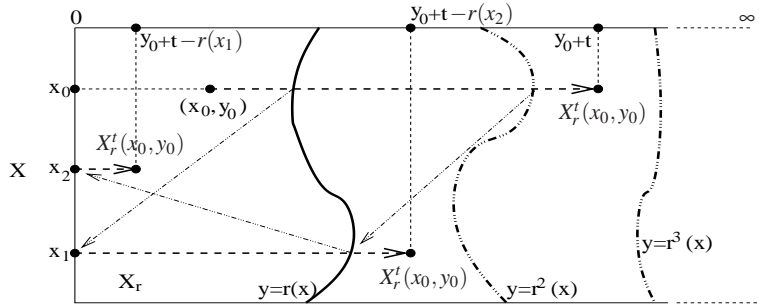}
  \caption{The equivalence relation defining the suspension
    flow of $f$ over the roof function $r$.}
  \label{fig:suspension}
\end{figure}

Set $M_r=\{ (x,y)\in M\times[0,+\infty): 0\le y < r(x) \}$
and $X^0$ the identity on $M_r$, where $M$ is a compact
Riemannian manifold.  For $x=x_0\in M$ denote by $x_n$ the
$n$th iterate $f^n(x_0)$ for $n\ge0$.  Denote $S_n^f \vfi(x_0)
= \sum_{j=0}^{n-1} \vfi( x_j )$ for $n\ge1$ and for any
given real function $\vfi$.  Then for each pair
$(x_0,s_0)\in X_r$ and $t>0$ there exists a unique $n\ge1$
such that $S_n r(x_0) \le s_0+ t < S_{n+1} r(x_0)$ and
define (see Figure~\ref{fig:suspension})
\begin{align*}
   X^t(x_0,s_0) = \big(x_n,s_0+t-S_n r(x_0)\big).
\end{align*}

The study of suspension (or special) flows is motivated by
modeling a flow admitting a cross-section. Such flow is
equivalent to a suspension semiflow over the Poincar\'e
return map to the cross-section with roof function given by
the return time function on the cross-section. This is a
main tool in the ergodic theory of uniformly hyperbolic
flows developed by Bowen and Ruelle~\cite{BR75}.

\subsubsection{Conditions on the base dynamics}
\label{sec:condit-base-dynamics}

We assume that the singular set $\cS$ (containing the points
where $f$ is either \emph{not defined}, \emph{discontinuous}
or \emph{not differentiable}) is regular, e.g. a submanifold
of $M$, and that $f$ is \emph{non-uniformly expanding}:
there exists $c>0$ such that for Lebesgue almost every $x\in
M$
\begin{align*}
  \limsup_{n\to+\infty}\frac1n S_n\psi(x) \le -c
\quad\mbox{where}\quad
\psi(x)=\log\big\| Df(x)^{-1} \big\|.
\end{align*}

Moreover we assume that $f$ has \emph{exponentially slow
  recurrence to the singular set $\cS$} i.e. for all
$\epsilon>0$ there is $\delta>0$ s.t.
\begin{align*}
  \limsup_{n\to+\infty}\frac1n\log\Leb\left\{x\in M:
    \frac1nS_n\big| \log d_{\delta}(x,\cS) \big|
    >\epsilon\right\}<0,
\end{align*}
where $d_\delta(x,y)=\dist(x,y)$ if $\dist(x,y)<\delta$ and
$d_\delta(x,y)=1$ otherwise.

These conditions ensure~\cite{ABV00} in particular the
existence of finitely many ergodic absolutely continuous (in
particular \emph{physical}) $f$-invariant probability
measures $\mu_1,\dots,\mu_k$ whose basins cover the manifold
Lebesgue almost everywhere.

We say that an $f$-invariant measure $\mu$ is an
\emph{equilibrium state} with respect to the potential $\log
J$, where $J=|\det Df|$, if $h_\mu(f)=\mu(\log J)$, that is
if \emph{$\mu$ satisfies the Entropy Formula}. Denote by
$\EE$ the family of all such equilibrium states.
It is not difficult to see that each physical measure in our
setting belongs to $\EE$.

We assume that \emph{ $\EE$ is formed by a unique
    absolutely continuous probability measure}.

\subsubsection{Conditions on the roof function}

We assume that $r:M\setminus\cS\to\RR^+$ has
\emph{logarithmic growth near $\cS$}: there exists
$K=K(\vfi)>0$ such that\footnote{$B(\cS,\delta)$ is the
  $\delta$-neighborhood of $\cS$.} $r\cdot\chi_{B(\cS,\delta)}\le
K\cdot\big| \log d_{\delta}(x,\cS) \big|$ for all small
enough $\delta>0$. We also assume that $r$ is bounded from
below by some $r_0>0$.

Now we can state the result on large deviations.

\begin{theorem}
  Let $X^t$ be a suspension semiflow over a non-uniformly
  expanding transformation $f$ on the base $M$, with roof
  function $r$, satisfying all the previouly stated conditions.

  Let $\psi:M_r\to\RR$ be continuous, $\nu=\mu\ltimes\Leb^1$
  be the induced invariant measure\footnote{for any
    $A\subset M_r$ set $\nu(A)=\mu(r)^{-1}\int d\mu(x)
    \int_0^{r(x)}\! ds \,\chi_A(x,s).$ } for the
  semiflow $X^t$ and $\lambda=\Leb\ltimes\Leb^1$ be the
  natural extension of volume to the space $M_r$.  Then
\begin{align*}
\limsup_{T\to\infty}\frac1T\log\lambda\left\{
z\in M_r: \left|
\frac1T\int_0^T \psi\left(X^t(z)\right)\, dt - \nu(\psi)
\right|>\epsilon 
\right\}<0.
\end{align*}
\end{theorem}

\subsubsection{Consequences for the Lorenz flow}

Now consider a Lorenz geometric flow as constructed in
Section \ref{sec:geoLorenz} and let $f$ be the
one-dimensional map associated, obtained quotienting over
the leaves of the stable foliation, see Figure
\ref{fig:GeoLorenz}.  This map has all the properties stated
previously for the base transformation. The Poincar\'e
return time gives also a roof function with logarithmic
growth near the singularity line.

The uniform contraction along the stable leaves implies that
the \emph{time averages of two orbits on the same stable
  leaf under the first return map are uniformly close} for
all big enough iterates. 
If $P:S\to[-1,1]$ is the
projection along stable leaves
\begin{lemma}
  For $\vfi:U\supset\Lambda\to\RR$ continuous and bounded,
  $\epsilon>0$ and $\vfi(x)=\int_0^{r(x)}\psi(x,t)\,dt$,
  there exists $\zeta:[-1,1]\setminus\cS\to\RR$ with
  logarithmic growth near $\cS$ such that
$
\Big\{\big|\frac1n
S_n^{R}\vfi-\mu(\vfi)\big|>2\epsilon\Big\}$
is contained in 
$$
  P^{-1}\Big(
  \big\{\big|\frac1n S^f_n\zeta-\mu(\zeta)\big|>\epsilon\big\}
  \cup
  \big\{
  \frac1n S^f_n\big|\log\dist_\delta(y,\cS)\big|> \epsilon
  \big\}
\Big).
$$
\end{lemma}
Hence in this setting it is enough to study the quotient map
$f$ to get information about deviations for the Poincar\'e return map. 
Coupled with the main result we are then able to
deduce

\begin{corollary}
  Let $X^t$ be a flow on $\RR^3$ exhibiting a Lorenz or a
  geometric Lorenz attractor with trapping region
  $U$. Denoting by $\Leb$ the normalized restriction of the
  Lebesgue volume measure to $U$, $\psi:U\to\RR$ a
  bounded continuous function and $\mu$ the unique physical
  measure for the attractor, then for any given $\epsilon>0$
  \begin{align*}
      \limsup_{T\to\infty}\frac1T\log\Leb\left\{
        z\in U: \left|
          \frac1T\int_0^T \psi\left(X^t(z)\right)\, dt - \mu(\psi)
        \right|>\epsilon 
      \right\}<0.
  \end{align*}
  Moreover for any compact $K\subset U$ such that $\mu(K)<1$
  we have
\begin{align*}
\limsup_{T\to+\infty}\frac1T \log \Leb\Big( \left\{ x\in K :
  X^t(x)\in K, 0<t<T \right\} \Big) < 0.
\end{align*}
\end{corollary}

\subsubsection{Idea of the proof}

We use properties of non-uniformly expanding
transformations, especially a large deviation bound recently
obtained~\cite{araujo-pacifico2006}, to deduce a large
deviation bound for the suspension semiflow reducing the
estimate of the volume of the deviation set to the volume of
a certain deviation set for the base transformation.

The initial step of the reduction is as follows.
For a continuous and bounded $\psi:M_r\to\RR$,  $T>0$
and $z=(x,s)$ with $x\in M$ and $0\le s < r(x) <\infty$,
 there exists the \textbf{lap number}
$n=n(x,s,T)\in\NN$ such that $S_{n}r(x)\le s+T < S_{n+1}
r(x)$,  and we can write
\begin{align*}
\int_0^T
\hspace{-0.2cm}\psi\big(X^t(z)\big)\,dt
=
\int_s^{r(x)}
\hspace{-0.6cm}\psi\big(X^t(x,0)\big)\,dt
&+\int_0^{T+s-S_{n}r(x)}
\hspace{-0.6cm}\psi\big(X^t(f^n(x),0)\big)\,dt
\\
&+\sum_{j=1}^{n-1} \int_0^{r(f^j(x))}
\hspace{-0.6cm}\psi\big(X^t(f^j(x),0)\big)\,dt.
\end{align*}
Setting $\vfi(x)=\int_0^{r(x)}\psi(x,0)\,dt$ we can rewrite
the last summation above as $S_n\vfi(x)$.  We get the
following expression for the time average
\begin{align*}
  \frac1T\int_0^T \hspace{-0.2cm}\psi\big(X^t(z)\big)\,dt
  =\frac1T S_n\vfi(x) 
  &
  -\frac1T\int_0^s\psi\big(X^t(x,0)\big)\,dt
  \\
  &+
  \frac1T\int_0^{T+s-S_n r(x)}
  \hspace{-0.6cm}\psi\big(X^t(f^n(x),0)\big)\,dt.
\end{align*}
Writing $I=I(x,s,T)$ for the sum of the last two integral
terms above, observe that for $\omega>0$, $0\le s<r(x)$ and
$n=n(x,s,T)$
\begin{align*}
  \left\{ (x,s)\in M_r : \left| \frac1T S_n\vfi(x) + I(x,s,T)
    -\frac{\mu(\vfi)}{\mu(r)}\right| > \omega\right\}
\end{align*}
is contained in
\begin{align*}
  \left\{ (x,s)\in M_r : \left| \frac1T S_n\vfi(x)
      -\frac{\mu(\vfi)}{\mu(r)}\right| > \frac\omega2 \right\}
\cup
\left\{(x,s)\in M_r: I(x,s,T) >\frac\omega2\right\}.
\end{align*}
The left hand side above is a \emph{deviation set for the
  observable $\vfi$ over the base transformation}, while the
right hand side will be \emph{bounded by the geometric
  conditions on $\cS$} and by a \emph{deviations bound for
  the observable $r$ over the base transformation}.

Analysing each set using the conditions on $f$ and $r$ and
noting that for $\mu$- and $\Leb$-almost every $x\in M$ and every
$0\le s<r(x)$
\begin{align*}
 \frac{S_nr(x)}n\le
 \frac{T+s}{n}\le\frac{S_{n+1}r(x)}n
 \quad\text{so}\quad
 \frac{n(x,s,T)}T\xrightarrow[T\to\infty]{}\frac1{\mu(r)},
\end{align*}
we are able to obtain the asymptotic bound of the Main
Theorem.

Full details of the proof are presented in~\cite{araujo2006a}.

The interested reader can find the proofs of the results
mentioned above in the papers listed below, the references
therein, and also in one of IMPA's texts~\cite{AraPac07} for
the XXV Brazilian Mathematical Colloquium.

\subsection*{Acknowledgments}
V.A. was partially supported by FCT-CMUP(Portugal) and CNPq,
FAPERJ and PRONEX-Dyn.Systems(Brazil). M.J.P. was partially
supported by CNPq, FAPERJ and PRONEX-Dyn.Systems(Brazil).
This work was done while M. J. Pacifico was enjoying a
post-doc leave at CRM Ennio De Giorgi-Scuola Normale
Superiore di Pisa with a CNPq fellowship and also partially
supported by CRM Ennio De Giorgi.
 

\def\cprime{$'$}

 \bibliographystyle{plain}

\end{document}